\documentclass{amsart}
\usepackage{amsmath}
\usepackage{amsfonts}
\usepackage{amssymb}

\newcommand{\QQ}{\mathbf{Q}}
\newcommand{\ZZ}{\mathbf{Z}}
\newcommand{\FF}{\mathbf{F}}
\newcommand{\EE}{\mathbf{E}}
\newcommand{\OO}{\widetilde{O}}
\newcommand{\ep}{\varepsilon}
\DeclareMathOperator{\LCM}{LCM}

\newcommand{\norm}[1]{{\hat{#1}}}   

\newtheorem{thm}{Theorem}
\newtheorem{lem}[thm]{Lemma}
\newtheorem{prop}[thm]{Proposition}
\theoremstyle{definition}
\newtheorem{example}[thm]{Example}
\theoremstyle{remark}
\newtheorem*{rem}{Remark}

\begin{document}

\title{Efficient computation of $p$-adic heights}
\author{David Harvey}
\begin{abstract}
We analyse and drastically improve the running time of the algorithm of Mazur, Stein and Tate for computing the canonical cyclotomic $p$-adic height of a point on an elliptic curve $E/\QQ$, where $E$ has good ordinary reduction at $p \geq 5$.
\end{abstract}

\maketitle


\section{Introduction}
\label{sec:intro}

Mazur, Stein and Tate \cite{MST} recently introduced a fast algorithm for computing the canonical, cyclotomic $p$-adic height $h_p(P)$ of a rational point $P$ on an elliptic curve $E/\QQ$ --- in this paper, referred to simply as the $p$-adic height --- in the case that $E$ has good ordinary reduction at $p \geq 5$. Their algorithm finds applications in, for example, numerical investigation of phenomena related to the $p$-adic Birch and Swinnerton-Dyer conjectures, since the $p$-adic height is part of the definition of the $p$-adic regulator term.

\cite{MST} did not give a bound for the running time of their algorithm; rather, the point of their paper is that the $p$-adic height can be computed feasibly at all, to reasonably high $p$-adic precision. In this paper we sketch an estimate for the running time of the algorithm of \cite{MST}, and we give several improvements that drastically improve the asymptotic running time, both for large $p$ and for high precision. With a contemporary desktop computer, it becomes straightforward to calculate a handful of $p$-adic digits of the height when $p$ is around $10^{11}$; and in the other direction, for small primes, one easily obtains thousands of $p$-adic digits.

We also carefully analyse the amount of $p$-adic precision that must be maintained throughout the calculation to guarantee a given number of correct digits of output; \cite{MST} did not discuss this issue in much detail. Obtaining a sharp bound is particularly important when $p$ is large, since computing additional unnecessary digits in intermediate calculations becomes very expensive.

In our running time estimates, we will generally ignore logarithmic factors, expressing all estimates in terms of the `soft-oh' notation $\OO(X)$, which means $O(X (\log X)^k)$ for some integer $k \geq 0$.

To state the main results, we introduce some notation. The elliptic curve $E/\QQ$ of interest is given by the Weierstrass equation
\begin{equation} \label{eq:weierstrass}
 y^2 + a_1 xy + a_3 y = x^3 + a_2 x^2 + a_4 x + a_6,
\end{equation}
with the usual invariant differential $\omega = dx/(2y + a_1 x + a_3)$. We assume that the $a_k$ are integers, but we do not assume that the equation is minimal. We assume that $p \geq 5$ and that $E$ has good ordinary reduction at $p$. We treat $E$ as fixed for the purposes of estimating running times.

The following result is proved in \S\ref{sec:E2}.
\begin{thm}
\label{thm:E2}
Let $N \geq 1$. The value of the $p$-adic modular form $\EE_2(E, \omega)$ may be computed mod $p^N$ in time $\OO(p N^2)$. If $p > 6N$, it may be computed in time $\OO(p^{1/2} N^{5/2})$.
\end{thm}

The next result, proved in \S\ref{sec:sigma}, concerns the $p$-adic sigma function $\sigma_p(t) \in \ZZ_p[[t]]$, as defined in \cite{padic-sigma}. This is a power series
 \[ \sigma_p(t) = t + c_2 t^2 + c_3 t^3 + \cdots, \]
where $t = -x/y$ is a parameter for the formal group of $E$. For any $r \geq 1$, let $I_r$ be the ideal of $\ZZ_p[[t]]$ given by
 \[ I_r = (p^r, p^{r-1} t, \ldots, p t^{r-1}, t^r). \]

\begin{thm}
\label{thm:sigma}
Let $N \geq 4$. Given the value of $\EE_2(E, \omega)$ mod $p^{N-3}$, the $p$-adic sigma function may be computed mod $I_N$ in time $\OO(N^2 \log p)$.
\end{thm}
Computing $\sigma_p(t)$ mod $I_N$ simply means that we are computing the coefficient of $t^k$ modulo $p^{N-k}$, for each $k = 1, \ldots, N-1$. The running time is optimal up to logarithmic factors, since the amount of data needed to represent $\sigma_p(t)$ mod $I_{N}$ is proportional to $N^2 \log p$. Note that for $N \leq 3$, computing $\sigma_p(t)$ mod $I_N$ is trivial, since $c_2$ is given simply by $a_1/2$ (see \S\ref{sec:sigma}).

Finally, we consider the problem of computing the $p$-adic height $h_p(P)$ mod $p^M$ of a point $P \in E(\QQ)$ for some $M \geq 2$, given $\sigma_p(t)$ as input. Following a suggestion of Christian Wuthrich, in this paper we normalise the $p$-adic height differently from \cite{MST}; our height is equal to the \cite{MST} height multiplied by the factor $2p$.

Let $n_1 = \# E(\FF_p)$, let $n_2$ be the least common multiple of the Tamagawa numbers of $E$, and let $n = \LCM(n_1, n_2)$. Put
 \[ M' = M + 2 v_p(n), \]
where $v_p$ denotes the usual additive $p$-adic valuation. Note that $v_p(n_1) \leq 1$, and we are assuming that $n_2 = O(1)$ since $E$ is fixed, so $M' = M + O(1)$.

The following result is proved in \S\ref{sec:height}. We denote by $C_P$ the amount of data required to represent the coordinates of $C_P$.

\begin{thm}
\label{thm:height}
Let $P \in E(\QQ)$, and suppose that $\sigma_p(t)$ is known mod $I_{M'+1}$. Then $h_p(P)$ may be computed mod $p^M$ in time $\OO(C_P + M \log^2 p + M^2 \log p)$.
\end{thm}

\subsection*{Organisation of the paper}

In \S\ref{sec:mst} we outline the Mazur--Stein--Tate algorithm, and indicate the steps in their algorithm that we intend to accelerate. In \S\ref{sec:E2}, \S\ref{sec:sigma} and \S\ref{sec:height} we prove Theorems \ref{thm:E2}, \ref{thm:sigma} and \ref{thm:height} respectively, and give several detailed examples of the various components of the algorithm. In \S\ref{sec:example} we give a higher-level example of the whole algorithm in operation. Finally in \S\ref{sec:samples} we give some examples of timings for an implementation of the algorithm.

\subsection*{Acknowledgements}

I would like to thank first of all William Stein for introducing me to the problem of computing $p$-adic heights, and for making many helpful suggestions on an early version of this paper. At the MSRI workshop ``Computing with Modular Forms'' (August 2006), I benefited greatly from working with Robert Bradshaw, Jennifer Balakrishnan and Liang Xiao on an implementation of the original Mazur--Stein--Tate algorithm in the computer algebra system SAGE \cite{sage}. A discussion with Christian Wuthrich regarding portions of his thesis helped cement the ideas for \S\ref{sec:division-polys}. Thanks also to Kiran Kedlaya for pointing out the method described in \cite{bernstein} for computing $p$-adic logarithms, and for many interesting conversations about his point-counting algorithm.


\section{A sketch of the Mazur--Stein--Tate algorithm}
\label{sec:mst}

In this section we sketch the original algorithm of \cite{MST}, and indicate the bottlenecks in the time complexity that this paper seeks to address.

The key insight of \cite{MST} is that it is possible to compute $\EE_2(E, \omega)$ efficiently using Kedlaya's algorithm \cite{kedlaya}, and then to use the value of $\EE_2(E, \omega)$ to deduce the $p$-adic sigma function. The $p$-adic height of a point $P \in E(\QQ)$ is then given in terms of $\sigma_p(t)$ by a simple formula. Their method for computing $\EE_2(E, \omega)$ is discussed in the proof of Theorem \ref{thm:E2} in \S\ref{sec:E2}. We now consider the other steps.

\subsection{Computing the $p$-adic sigma function}

The algorithm uses the fact that $\sigma_p(t)$ is the unique odd function in $\ZZ_p[[t]]$ of the form $\sigma_p(t) = t + \cdots$ satisfying the differential equation
\begin{equation}
\label{eq:sigma-de}
 x(t) + c = -\frac{d}{\omega}\left( \frac1{\sigma} \frac{d\sigma}{\omega} \right),
\end{equation}
where $c$ is the constant
 \[ c = \frac{a_1^2 + 4a_2 - \EE_2(E, \omega)}{12}, \]
and where $x(t) \in \ZZ_p[[t]]$ is the power series expansion of $x$ at the origin. Since $c$ is known from $\EE_2(E, \omega)$, it becomes a matter of computing $x(t)$, and then solving \eqref{eq:sigma-de} for $\sigma_p(t)$.

The series $x(t)$ is determined from an auxiliary power series $w(t) = \sum_{n \geq 0} s_n t^n$, using a certain recursive formula to compute the $s_n$. A bottleneck arises here, in that the recursive formula for $s_n$ involves a term of the form $\sum_{i+j+k = n} s_i s_j s_k$, which requires $O(n^2)$ ring operations to evaluate. To compute the $p$-adic height to precision $p^N$, it turns out to be necessary to compute $O(N)$ coefficients, to precision about $p^N$ each. Therefore computing $x(t)$ to the desired precision requires $O(N^3)$ ring operations, for a time cost proportional to at least $N^3 \log(p^N) = N^4 \log p$. In \S\ref{sec:sigma} we will give a different method for computing $x(t)$ whose running time is only $\OO(N^2 \log p)$.

To solve \eqref{eq:sigma-de}, they first use the formal logarithm to to change variables from $t$ to the parameter $z$ on the \emph{additive} group. After making this substitution, the differential equation takes the simpler form
 \[ x(z) + c = -\frac{d}{dz} \left( \frac1{\sigma(z)} \frac{d\sigma}{dz} \right), \]
which can be solved by integration and exponentation of power series. A bottleneck arises here also: to perform the change of variables, it is necessary to invoke several power series composition and reversion operations. To the author's knowledge, the Brent--Kung algorithm \cite{brent-kung} is the best algorithm available for computing series reversions and compositions. It has complexity $\OO(n^{3/2})$ in the number of terms, so when working to precision $p^N$ the running time would be proportional to at least $N^{5/2} \log p$. In \S\ref{sec:sigma} we show that it is entirely unnecessary to change variables; the original equation \eqref{eq:sigma-de} can be solved directly in time $\OO(N^2 \log p)$.

\subsection{Computing the $p$-adic height}
\label{sec:mst-2}

Let $P \in E(\QQ)$ be a nonzero point. Its affine coordinates may be written uniquely in the form
 \[ P = (x(P), y(P)) = \left(\frac{\alpha(P)}{d(P)^2}, \frac{\beta(P)}{d(P)^3}\right), \]
where $(\alpha(P), d(P)) = (\beta(P), d(P)) = 1$ and $d(P) \geq 1$.

First \cite{MST} consider the case that $P$ satisfies two conditions:
\begin{itemize}
\item[(A1)] $P$ reduces to the zero element of $E(\FF_p)$; and
\item[(A2)] $P$ reduces to a nonsingular point of $E(\FF_\ell)$ for all primes $\ell$ at which $E$ has bad reduction.
\end{itemize}
Condition (A1) implies that $t(P) = -x(P)/y(P)$ is divisible by $p$. Under these conditions, $h_p(P)$ is given by the formula
\begin{equation}
\label{eq:height-formula}
 h_p(P) = 2 \log_p\left(\frac{\sigma(P)}{d(P)} \right),
\end{equation}
where $\log_p$ is the Iwasawa $p$-adic logarithm, and where $\sigma(P) = \sigma(t(P))$. (Recall that our normalisation for $h_p(P)$ differs from that of \cite{MST} by a factor of $2p$.)

To handle arbitrary nonzero $Q \in E(\QQ)$, they proceed as follows. As in Theorem \ref{thm:height}, let $n_1 = \# E(\FF_p)$, let $n_2$ be the least common multiple of the Tamagawa numbers of $E$, and let $n = \LCM(n_1, n_2)$. Then $P = nQ$ \emph{does} satisfy (A1) and (A2), so one may use \eqref{eq:height-formula} to compute $h_p(P)$. From this one obtains $h_p(Q) = h_p(P)/n^2$, since $h_p$ is a quadratic form.

The above procedure has a serious bottleneck when $p$ is large. The difficulty is that one must actually compute the coordinates of $nQ$, which will generally have about $n^2$ digits (assuming that $Q$ is non-torsion). From the Weil bound we have $\# E(\FF_p) \approx p$, so $n$ is roughly proportional to $p$. Therefore the time complexity is at least proportional to $p^2$, and for sufficiently large $p$ will dwarf even the time required to compute $\EE_2(E, \omega)$. In \S\ref{sec:height} we will show how to avoid this problem by judicious use of the division polynomials associated to $E$, achieving a running time only polylogarithmic in $p$.


\section{Computing $\EE_2(E, \omega)$}
\label{sec:E2}

In this section we prove Theorem \ref{thm:E2}: given $N \geq 1$, we wish to compute $\EE_2(E, \omega)$ to precision $p^N$.

The first step is to compute the matrix $F$ of the absolute $p$-th power Frobenius on the basis $\{dx/y, x\, dx/y\}$ of the Monsky--Washnitzer cohomology of $E$, to precision $p^N$. This may be done using either Kedlaya's algorithm \cite{kedlaya} in time $\OO(p N^2)$, or if $p > 6N$ one may use a modification of Kedlaya's algorithm \cite{kedlaya-large-p} that has running time $\OO(p^{1/2} N^{5/2})$. The optimal crossover point between the two algorithms would depend on the implementations; one expects it to be roughly of the form $p = CN$ for some constant $C$.

As explained in \cite[\S 3.2]{MST}, the value of $\EE_2(E, \omega)$ may be deduced from the matrix as follows. Write
 \[ F^N = \begin{pmatrix} A & B \\ C & D \end{pmatrix}. \]
Then
 \[ \EE_2(E, \omega) = -12 B/D \pmod{p^N}. \]
Note that $D$ is a unit (see \cite[\S 3.2]{MST}), so the result is correct mod $p^N$. Using repeated squaring, computing $F^N$ requires $O(\log N)$ matrix multiplications, each taking time $\OO(N \log p)$, so the time for this step is only $\OO(N \log p)$. 

\subsection{A constant factor improvement}
\label{sec:E2-improvement}

Kedlaya's algorithm was originally designed for hyperelliptic curves, but because we are in the special setting of an elliptic curve, it is possible to obtain a constant factor speedup. Indeed, the \emph{trace} of the matrix may be determined by any of various fast algorithms for counting $\#E(\FF_p)$, and the determinant is known to be $p$. Since Kedlaya's algorithm computes the two columns of the matrix independently, the idea is simply to compute only one column, and fill in the other column from knowledge of the trace and determinant. (Alternatively, one could use the trace and determinant as a strong correctness check.) This idea was communicated to the author by William Stein, who attributes it to Eyal Goren.

The speedup will be at most a factor of two. In practice it will be less, since it only affects ``Step 2'' of Kedlaya's algorithm, not ``Step 1'' (in the language of \cite[\S4]{kedlaya}). It only applies to Kedlaya's original algorithm; it does not apply to the algorithm of \cite{kedlaya-large-p}.

In practice some care is needed to avoid precision loss. A particularly badly conditioned example is given by the curve $y^2 = x^3 + 7x + 8$ at the prime $p = 11$, whose Frobenius matrix is given to precision $O(11^3)$ by
 \[  F = \begin{pmatrix}
          11 \cdot 104 + O(11^3)  &  11 \cdot 16 + O(11^3) \\
          11^2 \cdot 7 + O(11^3)  &  185 + O(11^3)  \end{pmatrix}.  \]
If the trace and determinant are known, then the second column mod $11^3$ only determines the first column mod $11^2$, and the first column mod $11^3$ only determines the second column mod $11$.

This problem may be avoided by a simple change-of-basis argument, as follows. One first uses Kedlaya's algorithm to compute the second column only mod $p$. The running time for this step is $\OO(p)$. If the top-right entry is a unit, then no precision loss will arise; one simply uses Kedlaya's algorithm at full precision $p^N$ to compute the second column, and this suffices to determine the first column mod $p^N$ from the trace and determinant.

In the unlucky event that the top-right entry is not a unit, one instead runs Kedlaya's algorithm at precision $p^N$ with respect to the basis $\{dx/y, (1+x)dx/y\}$. In other words, one applies Kedlaya's reduction algorithm to the image under Frobenius of the differential $(1+x) dx/y$, and expresses the result as a linear combination of $dx/y$ and $(1+x)dx/y$. The top-right entry of the matrix on this new basis is then guaranteed to be a unit, so the second column is sufficient to determine the first column mod $p^N$ (and the trace does not depend on the choice of basis). Conjugating by the change-of-basis matrix then yields the matrix on the original basis, to full precision.


\section{Computing the $p$-adic sigma function}
\label{sec:sigma}

This section is devoted to proving Theorem \ref{thm:sigma}. We assume that $N \geq 4$ and that $c = (a_1^2 + 4a_2 - \EE_2)/12$ is known mod $p^{N-3}$. We wish to compute $\sigma_p(t)$ mod $I_N$, in time $\OO(N^2 \log p)$.

Throughout this section we assume that asymptotically fast polynomial arithmetic is being used; specifically, that multiplication and division of polynomials of degree $d$ over $\ZZ/p^r\ZZ$ may be accomplished in time $\OO(dr\log p)$.

Let $x(t)$ and $y(t)$ be the power series expansions of $x$ and $y$ around the origin. Let $w(t) = -1/y(t)$, and let
 \[ \omega(t) = \frac{x'(t)}{2y(t) + a_1 x(t) + a_3}, \]  
so that $\omega(t) dt$ is the series expansion of the invariant differential $\omega$. The first few terms of each series are given by 
\[ \begin{aligned}
 x(t) & = t^{-2} - a_1 t^{-1} - a_2 + \cdots, \\
 y(t) & = -t^{-3} + a_1 t^{-2} + a_2 t^{-1} + \cdots, \\
 w(t) & = t^3 + a_1 t^4 + (a_1^2 + a_2) t^5 + \cdots, \\
 \omega(t) & = 1 + a_1 t + (a_1^2 + a_2) t^2 + \cdots.
\end{aligned} \]
(See also \cite[Ch.~IV]{silverman}, which discusses these expansions in some detail, using slightly different notation.)

Let $R = \ZZ/p^{N-3}\ZZ$. Our first task is to compute the above series over $R$, with $N+1$ terms each.

\begin{prop}
\label{prop:newton}
The series $t^2 x(t)$, $t^3 y(t)$, $t^{-3} w(t)$ and $\omega(t)$ may be computed up to $O(t^{N+1})$, with coefficients in $R$, in time $\OO(N^2 \log p)$.
\end{prop}
\begin{proof}
The series $w(t)$ satisfies the algebraic equation
 \[ w = t^3 + a_1 tw + a_2 t^2 w + a_3 w^2 + a_4 tw^2 + a_6 w^3, \]
and so it may be solved using Newton's method. We start with the initial approximation $w(t) = t^3$, and then repeatedly apply the Newton iteration
\begin{equation}
\label{eq:w-newton-iteration}
\begin{aligned}
 w' & = w - \frac{w - t^3 - a_1 tw - a_2 t^2 w - a_3 w^2 - a_4 tw^2 - a_6 w^3}
                 {1 - a_1 t - a_2 t^2 - 2a_3 w - 2a_4 tw - 3a_6 w^2} \\
      & = \frac{t^3 - a_3 w^2 - a_4 t w^2 - 2a_6 w^3}
          {1 - a_1 t - a_2 t^2 - 2a_3 w - 2a_4 tw - 3a_6 w^2}.
\end{aligned}
\end{equation}
The arithmetic is performed using truncated power series in $R[t]$. It is straightforward to check that the number of correct terms doubles with each iteration. Each iteration requires several power series multiplications and one inversion. As is typical in applications of Newton's method to power series, the total time to obtain $n$ terms is a constant multiple of the time required for a length $n$ polynomial multiplication. Therefore the total time to compute $O(N)$ terms of $w(t)$ with coefficients in $R$ is $\OO(N \log(p^{N-3})) = \OO(N^2 \log p)$.

Given $t^{-3} w(t)$ up to $O(t^{N+1})$, we may then deduce $x(t) = t/w(t)$, $y(t) = -1/w(t)$ and $\omega(t) = x'(t)/(2y(t) + a_1 x(t) + a_3)$, also with coefficients in $R$, to the desired number of terms, in time $\OO(N^2 \log p)$.
\end{proof}

Next we consider the differential equation \eqref{eq:sigma-de}, which may be rewritten as
 \[ x(t) + c = \frac{-1}{\omega(t)} \left( \frac{\sigma_p'(t)}{\sigma_p(t) \omega(t)} \right)'. \]
It is convenient to rephrase this in terms of the unit power series
 \[ \theta(t) = t^{-1} \sigma_p(t) = 1 + \cdots \in \ZZ_p[[t]]. \]
To solve for $\sigma_p(t)$ mod $I_N$ it suffices to solve for $\theta(t)$ mod $I_{N-1}$. We have $\theta'(t)/\theta(t) = \sigma_p'(t)/\sigma_p(t) + t^{-1}$, so $\theta(t)$ satisfies the equation
\begin{equation}
\label{eq:theta-2}
 x(t) + c = \frac{-1}{\omega(t)} \left(\frac1{\omega(t)} \left(\frac1t + \frac{\theta'(t)}{\theta(t) }\right) \right)'.
\end{equation}
Manipulating this equation formally, we obtain
\begin{equation}
\label{eq:theta}
 \frac{\theta'(t)}{\theta(t)} = h(t) \qquad \in \ZZ_p[[t]],
\end{equation}
where
\begin{equation}
\label{eq:H-formula}
 h(t) =  -\frac1t - \omega(t) \left( \int (x(t) + c) \omega(t)dt + C\right)
\end{equation}
for some constant of integration $C$. (The formal integration operator is assumed to output a series with zero constant term.)

The constant $C$ is determined by the condition that $\sigma_p(t)$ is an \emph{odd} function. This condition means that $\sigma_p(i(t)) = -\sigma_p(t)$, where $i(t) = -t - a_1 t^2 + \cdots$ is the formal inverse law \cite[IV \S1]{silverman}. It implies that the coefficient of $t^2$ in $\sigma_p(t)$ is $a_1/2$, so that the coefficient of $t$ in $\theta(t)$ is $a_1/2$. Substituting this into \eqref{eq:theta} and using the expansions for $\omega(t)$ and $x(t)$ given earlier, one finds that $C = a_1/2$.

\begin{prop}
\label{prop:h-approx}
The series $h(t)$ may be computed, with the coefficient of $t^j$ known mod $p^{N-3-\lfloor \log_p(j)\rfloor}$ for $1 \leq j < N$, and with the constant term known mod $p^{N-2}$, in time $\OO(N^2 \log p)$.
\end{prop}
In the above proposition, note that $\log_p(j)$ refers to the base $p$ logarithm, not the $p$-adic logarithm.
\begin{proof}
First we compute $(x(t) + c) \omega(t)$. From Proposition \ref{prop:newton} this series is known up to $O(t^{N-1})$ with coefficients in $R = \ZZ/p^{N-3}\ZZ$. To integrate it, we must check that the coefficient of $t^j$ has valuation at least $v_p(j+1)$ for each $j$; this follows for example from \eqref{eq:theta-2}. (This fact is the basis for the `integrality algorithm' of \cite{MST}.)

After integrating, the series
 \[ \int (x(t) + c) \omega(t) dt + \frac{a_1}2 \]
is known up to $O(t^N)$, but the coefficient of $t^j$ is known only mod $p^{N-3-v_p(j)}$, because of the divisions by powers of $p$ in the integration step. Then we must multiply by $-\omega(t)$ and subtract $1/t$, to obtain
 \[ h(t) = - \frac1t - \omega(t) \left(\int (x(t) + c) \omega(t) dt + \frac{a_1}2\right). \]
This series is also known up to $O(t^N)$; the effect of multiplying by $\omega(t)$ is to propagate the errors to higher order terms, so now the coefficient of $t^j$ is known only mod $p^{N-3-\lfloor \log_p(j)\rfloor}$.

The constant term must be determined mod $p^{N-2}$ separately; from the explicit series expansions given above, it is simply $a_1/2$.
\end{proof}

Finally we must solve \eqref{eq:theta} for $\theta(t)$. Formally the solution is given by
 \[ \theta(t) = \exp\left(\int h(t) dt\right). \]
Computationally this is problematic because both the integration and exponentiation steps introduce denominators, and also some $p$-adic precision is lost. The following result gives an efficient means to solve equations of the form $F'/F = f$, without any denominators appearing in intermediate steps, and with good control over precision loss.

\begin{prop}
\label{prop:brent}
Let $k \geq 1$ and $1 \leq n < p^{k/2}$, and let $S = \ZZ/p^k\ZZ$.  Let $f \in S[t]/(t^{n-1})$, and suppose that there exists $F \in S[t]/(t^n)$, with $F(0) = 1$, such that $F'/F = f$. Then $F$ may be determined mod $J$ in time $\OO(nk \log p)$, where $J$ is the ideal of $S[t]/(t^n)$ given by
 \[ J = (p^{k-1} t^p, p^{k-2} t^{p^2}, \ldots). \]
\end{prop}
\begin{rem}
Note that $J$ is the ideal that captures the types of $p$-adic error terms that occur in a power series integration. Namely, if $g \in S[t]/(t^{n-1})$, and if the coefficient of $t^j$ in $g$ is divisible by $j+1$ for each $j$, then $g$ has an integral in $S[t]/(t^n)$. In general it is not uniquely determined, but it is determined at least mod $J$.
\end{rem}
\begin{proof}
The algorithm we describe is essentially that of Brent \cite{brent}, with some additional analysis to track the $p$-adic error terms.

We begin with an initial approximation $F_0(t) = 1 \in S[t]/(t^n)$. We will refine it iteratively, obtaining a sequence $F_i \in S[t]/(t^n)$ such that $F_i - F \in J + (t^{2^i})$ for each $i \geq 0$. After $\lceil \log_2 n \rceil$ steps, we will have $F_i - F \in J$ as desired. Each step is dominated by a polynomial multiplication and a division of length $2^i$, so the total time is a constant multiple of the time required for a single polynomial multiplication of length $n$, which is $\OO(n \log(p^k)) = \OO(nk \log p)$.

Now we explain the iterative step. Suppose that $F_i - F \in J + (t^{2^i})$. Since $F$ is invertible, we have
 \[ F_i = (1+\ep)F \]
for some $\ep \in J + (t^{2^i})$, and so
 \[ \frac{F_i'}{F_i} - f = \frac{F'}{F} + \frac{\ep'}{1+\ep} - f = \frac{\ep'}{1+\ep}. \]

Morally speaking we would like to integrate $F_i'/F_i - f$ to obtain $\log(1+\ep)$, but the latter does not make sense in our ring. Instead, we write
 \[ \frac{\ep'}{1+\ep} = \ep' - \ep'\ep + \frac{\ep' \ep^2}{1+\ep}. \]
The hypothesis $n < p^{k/2}$ implies that $J^2 = 0$, so that $\ep^2 \in t^{2^i}J$. We also have $\ep' \in J + (t^{2^i - 1})$, so $\ep' \ep^2 \in (t^{2^{i+1} - 1})$. Consequently
 \[ \frac{F_i'}{F_i} - f \in \ep' - \ep'\ep + (t^{2^{i+1} - 1}). \]
Therefore $F_i'/F_i - f$ may be integrated at least up to $O(t^{2^{i+1}})$; that is, we may compute a $G \in S[t]/(t^n)$ such that
 \[ G \in \ep - \frac{\ep^2}2 + J + (t^{2^{i+1}}). \]
(The extra $J$ term is introduced by errors in the integration.) Since $\ep^2 \in J$ we have in fact
 \[ G \in \ep + J + (t^{2^{i+1}}). \]

Now it is straightforward to define $F_{i+1}$; we simply take
 \[  F_{i+1} = F_i(1 - G). \]
From the above estimates this satisfies
 \[ \begin{aligned}
    F_{i+1} & \in F(1+\ep)(1-\ep) + J + (t^{2^{i+1}}) \\
            & = F - \ep^2 F + J + (t^{2^{i+1}}) \\
            & = F + J + (t^{2^{i+1}})
   \end{aligned} \]
as desired.
\end{proof}

Now we may complete the proof of Theorem \ref{thm:sigma}. Let $\hat h(t)$ be the approximation to $h(t)$ produced by Proposition \ref{prop:h-approx}, treated as an element of $\ZZ_p[[t]]$, and let 
 \[ \hat \theta(t) = \exp\left(\int \hat h(t) \, dt\right) \in \QQ_p[[t]]. \]
Observe that
 \[ \frac{\hat\theta(t)}{\theta(t)} = \exp\left(\int \hat h(t) - h(t) \, dt \right). \]
By Proposition \ref{prop:h-approx}, for $1 \leq j < N$ the coefficient of $t^j$ in $\hat h - h$ has valuation at least $N - 3 - \lfloor\log_p j \rfloor$, so the coefficient of $t^j$ in $\int \hat h - h$ has valuation at least
 \[ N - 3 - \lfloor\log_p (j-1) \rfloor - \lfloor\log_p j \rfloor. \]
Therefore the coefficient of $t^j$ in $\exp(\int \hat h - h) - 1$ has valuation at least
 \[ N - 3 - \lfloor\log_p (j-1) \rfloor - \lfloor\log_p j \rfloor - \left\lfloor \frac{j}{p-1} \right\rfloor, \]
because $v_p(j!) \leq j/(p-1)$. Since $p \geq 5$ we have
 \[ 3 + \lfloor\log_p (j-1) \rfloor + \lfloor\log_p j \rfloor + \left\lfloor \frac{j}{p-1} \right\rfloor \leq j+1 \]
for all $j \geq 2$. (It suffices to estimate the left hand side for $p = 5$. For large enough $j$, the $j/4$ term dominates; one must also check a few small values of $j$ directly.)

This shows that the $t^j$ coefficients of $\hat\theta(t)$ and $\theta(t)$ agree mod $p^{N-j-1}$ for $2 \leq j < N$. This holds for $j = 1$ also, because of the extra condition on the constant term of $h(t)$ in Proposition \ref{prop:h-approx}. In other words, we have $\hat\theta - \theta \in I_{N-1}$. In particular, the coefficients of $\hat\theta$ are integral up to $O(t^{N-1})$, and the hypotheses of Proposition \ref{prop:brent} are satisfied with $F = \hat\theta$, $f = \hat h$, $k = N-2$ and $n = N-1$.

The output of Proposition \ref{prop:brent} is $\hat\theta$ mod $J$. This determines $\theta$ mod $I_{N-1}$, except for the constant term, which is known to be $a_1/2$.

\begin{example}
\label{ex:sigma}

We will walk through the steps of computing $\sigma_p(t)$ for the curve
 \[ E: y^2 + xy + y = x^3 - 460x - 3830, \]
which is `26a2' in Cremona's database, and the prime $p = 5$. We will take $N = 9$; that is, we will determine the coefficients of $t, t^2, \ldots, t^8$ mod $5^8, 5^7, \ldots, 5^1$ respectively.

First we must compute $w(t)$ up to $O(t^{13})$, with coefficients mod $5^6$, using Proposition \ref{prop:newton}. The successive approximations are
\[
\begin{split}
w(t) & = t^3 + O(t^4), \\
w(t) & = t^3 + t^4 + t^5 + 2t^6 + 15169t^7 + O(t^8), \\
w(t) & = t^3 + t^4 + t^5 + 2t^6 + 15169t^7 + 14252t^8 + 9048t^9 \\
 & \qquad + 9516t^{10} + 9477t^{11} + 14344t^{12} + O(t^{13}).
\end{split}
\]
From this we obtain the other series attached to $E$,
\begin{multline*}
y(t) = -1/w(t) = -t^{-3} + t^{-2} + 1 + 15166t + 15166t^2 \\
      + 11337t^3 + 6589t^4 + 9397t^5 + 8273t^6 + O(t^7),
\end{multline*}
\begin{multline*}
x(t) = t/w(t) = t^{-2} - t^{-1} - t + 459t^2 + 459t^3 \\
      + 4288t^4 + 9036t^5 + 6228t^6 + 7352t^7 + O(t^8),
\end{multline*}
\begin{multline*}
\omega(t) = x'(t) / (2y(t) + a_1 x(t) + a_3) = 1 + t + t^2 + 3t^3 + 14712t^4 + 12878t^5 \\
    + 14267t^6 + 1881t^7 + 4058t^8 + 2267t^9 + O(t^{10}).
\end{multline*}

Now we apply Proposition \ref{prop:h-approx}. Let us assume that $\EE_2 = 4303 + O(5^6)$ has been computed in advance, so we have
 \[ c = \frac{1^2 + 4 \cdot 0 - 4303}{12} = 7454 + O(5^6), \]
and then
\begin{multline*}
 (x(t) + c)\omega(t) = t^{-2} + 7454 + 7455t + 6996t^2 + 5820t^3 + \\
      13590t^4 + 11924t^5 + 15504t^6 + 1081t^7 + O(t^8).
\end{multline*}
Observe that the $t^{-1}$ term is zero, so the series may be integrated. Our choice of $c$ ensures that the $t^4$ term is divisible by $p = 5$, so that the integral has coefficients in $\ZZ_5$:
\begin{multline*}
 \int (x(t) + c)\omega(t) dt + \frac{a_1}{2} = -t^{-1} + 7813 + 7454t + 11540t^2 + \\
    2332t^3 + 1455t^4 + 2718t^5 + 12404t^6 + 4447t^7 + 13807t^8 + O(t^9).
\end{multline*}
Note that the $t^5$ coefficient is now only correct mod $5^5$, because we lost a digit during the integration (in fact the correct coefficient is $12093 \pmod{5^6}$). We obtain
\begin{multline*}
 \hat h(t) = \frac{-1}{t} - \omega(t) \left(\int (x(t) + c)\omega(t) dt + \frac{a_1}{2} \right) \\
   = 7813 + 359t + 4446t^2 + 1197t^3 + 14708t^4 \\
       + 6580t^5 + 6770t^6 + 1524t^7 + 2441t^8 + O(t^9).
\end{multline*}
The preceding multiplication by $\omega(t)$ caused the incorrect digit to wash through to the higher order terms, so now the $t^5$ through $t^8$ coefficients are only correct mod $5^5$. Also, the constant term $a_1/2$ is only correct mod $5^6$, but we need to increase its precision to $5^7$, so $\hat h(t)$ becomes
\begin{multline*}
 \hat h(t) = 39063 + 359t + 4446t^2 + 1197t^3 + 14708t^4 + 6580t^5 \\
   + 6770t^6 + 1524t^7 + 2441t^8 + O(t^9).
\end{multline*}

Now we run Brent's algorithm (Proposition \ref{prop:brent}) to solve for $\hat\theta(t)$, with the polynomial arithmetic performed mod $5^7$. The successive approximations are
\[
\begin{split}
 \hat\theta(t) & = 1 + O(t), \\
 \hat\theta(t) & = 1 + 39063t + O(t^2), \\
 \hat\theta(t) & = 1 + 39063t + 68539t^2 + 12965t^3 + O(t^4), \\
 \hat\theta(t) & = 1 + 39063t + 68539t^2 + 12965t^3 + 30804t^4 + 14720t^5 \\
         & \qquad + 10063t^6 + 25830t^7 + O(t^8). \\
\end{split}
\]
Finally $\sigma_p(t)$ is obtained as $t\hat\theta(t)$, but Theorem \ref{thm:sigma} only guarantees the result is correct mod $I_N$, so we have
\begin{multline*}
 \sigma_p(t) = t + (39063 + O(5^7))t^2 + (6039 + O(5^6))t^3 \\
        + (465 + O(5^5))t^4 + (179 + O(5^4))t^5 + (95 + O(5^3))t^6 \\
        + (13 + O(5^2))t^7 + (0 + O(5))t^8 + O(t^9).
\end{multline*}
\end{example}


\section{Computing the $p$-adic height}
\label{sec:height}

In this section we prove Theorem \ref{thm:height}. Recall that $n_1 = \# E(\FF_p)$, $n_2$ is the least common multiple of the Tamagawa numbers of $E$, and $n = \LCM(n_1, n_2)$. We are given a point $P \in E(\QQ)$ as input, whose coordinates require space $C_P$ to store. Given the sigma function mod $I_{M'+1}$ as input, where
 \[ M' = M + 2v_p(n), \]
we wish to compute $h_p(P)$ to precision $p^M$, in time $\OO(C_P + M \log^2 p + M^2 \log p)$.

As noted in \S\ref{sec:mst}, computing $nP$ directly is not feasible for large $p$. Therefore we take a more indirect approach. In what follows, the symbols $\alpha(P)$, $\beta(P)$, $d(P)$, $x(P)$, $y(P)$ and $t(P)$ are the same as those defined in \S\ref{sec:mst-2}.

First we compute $Q = n_2 P$, so that $Q$ satisfies the condition (A2) defined in \S\ref{sec:mst-2}.

Now let $m = n / n_2$. Because $mQ$ ($= nP$) satisfies both conditions (A1) and (A2), we may express $h_p(mQ)$ using \eqref{eq:height-formula}. Since $h_p$ is a quadratic form, we obtain
\begin{equation}
\label{eq:height-2}
 h_p(P) = \frac{2}{n^2} \log_p\left(\frac{\sigma_p(mQ)}{d(mQ)}\right).
\end{equation}
To determine $h_p(P)$ mod $p^M$, it suffices to compute
\begin{equation}
\label{eq:sigma-subst}
 \log_p\left(\frac{\sigma_p(mQ)}{d(mQ)}\right)
        =  \log_p \left(\frac{-\alpha(mQ)}{\beta(mQ)}
               \left( 1 + \sum_{k \geq 1} c_{k+1} t(mQ)^k \right) \right)
\end{equation}
mod $p^{M'}$, where we have expanded $\sigma_p(t)$ as
 \[ \sigma_p(t) = t + c_2 t^2 + c_3 t^3 + \cdots. \]

Note that $\alpha(mQ) / \beta(mQ)$ is a unit, since both $\alpha(mQ)$ and $\beta(mQ)$ are relatively prime to $d(mQ)$, which is divisible by $p$. Also $1 + \sum_{k \geq 1} c_{k+1} t(mQ)^k$ is a unit, again because $t(mQ) = -d(mQ)\alpha(mQ)/\beta(mQ)$ is divisible by $p$.

\begin{lem} \label{lem:padic-log}
Suppose that $u \in \ZZ_p^*$ is a unit and that $N \geq 1$ is an integer. To determine $\log_p u$ mod $p^N$ it suffices to know $u$ mod $p^N$.
\end{lem}
\begin{proof}
From $\log_p u = \frac{1}{p-1} \log_p(u^{p-1})$ we may assume that $u \equiv 1 \pmod p$. Since $\log_p(u + \ep) = \log_p u  + \log_p(1 + \ep/u)$, the result amounts to checking that if $v_p(\ep) \geq N$, then $v_p (\log_p(1+\ep)) \geq N$. Using the power series expansion of $\log_p(1+x)$, this follows from the elementary estimate $v_p(\ep^n/n) \geq v_p \ep$, that is, $v_p(n) \leq n-1$ for all $n \geq 1$.
\end{proof}

By Lemma \ref{lem:padic-log}, it suffices to compute
 \[ \frac{\alpha(mQ)}{\beta(mQ)}
 \qquad \text{and} \qquad
  1 + \sum_{k \geq 1} c_{k+1} t(mQ)^k \]
mod $p^{M'}$. By hypothesis we have available the sigma function mod $I_{M'+1}$, which means that we know $c_j$ mod $p^{M'+1-j}$ for $2 \leq j \leq M'$. Since $p$ divides $t(mQ)$, it therefore suffices to compute $\alpha(mQ)/\beta(mQ)$ and $t(mQ)$ mod $p^{M'}$. For this, we have the following result, which is proved in \S\ref{sec:division-polys} below.

\begin{prop}
\label{prop:division-polys}
Suppose that $Q \in E(\QQ)$ reduces to a non-singular point of $E(\FF_\ell)$ for every prime $\ell$ at which $E$ has bad reduction. Let $R \geq 1$ be odd, and $m \geq 2$. Given the values of $\alpha(Q)$, $\beta(Q)$ and $d(Q)$ mod $R$, one may compute $\pm\alpha(mQ)$, $\beta(mQ)$ and $\pm d(mQ)$ mod $R$ in time $\OO(\log R \log m)$.

(The $\pm$ symbols indicate that $\alpha(mQ)$ and $d(mQ)$ will be correct only up to sign, and that the signs will agree.)
\end{prop}

Proposition \ref{prop:division-polys} completely avoids the problem of the coordinates of $mQ$ growing out of control, since the dependence on $m$ is only logarithmic.

We apply Proposition \ref{prop:division-polys} with $R = p^{M'}$, to determine $\alpha(mQ)/\beta(mQ)$ and $t(mQ)$ mod $p^{M'}$. The sign ambiguity is irrelevant, since in $t(mQ) = -d(mQ)\alpha(mQ)/\beta(mQ)$ the signs cancel out, and the sign of $\alpha(mQ)/\beta(mQ)$ is not needed since the $p$-adic logarithm is insensitive to the sign of its input.

Now we analyse the running time. Recall that $n_2$ is assumed to be $O(1)$, so the time for computing $Q = n_2 P$ is $\OO(C_P)$. Applying Proposition \ref{prop:division-polys} costs $\OO(\log(p^{M'}) \log m)$. We have $M' = O(M)$ and $m = O(p)$, so this is $\OO(M \log^2 p)$. We must substitute $t(mQ)$ into \eqref{eq:sigma-subst}, which requires $O(M')$ ring operations in $\ZZ/p^{M'}\ZZ$, costing time $\OO(M' \log(p^{M'}) = \OO(M^2 \log p)$. Finally we must compute the $p$-adic logarithm in \eqref{eq:height-2}. Using the series expansion of $\log_p(1+x)$, this requires $O(M')$ ring operations, for a cost of $\OO(M^2 \log p)$; with more effort one can obtain $\OO(M \log p)$ \cite[\S 16]{bernstein}.

\subsection{Proof of Proposition \protect{\ref{prop:division-polys}}}
\label{sec:division-polys}

Our main tools in the proof are the division polynomials $\psi_m$ associated to our choice of Weierstrass equation for $E$, and a non-cancellation result of Wuthrich \cite[Prop.~IV.2]{wuthrich} that in effect controls the amount of cancellation that can occur while computing $mQ$ from $Q$. For further background on division polynomials, including proofs of the assertions we use in this section, we refer the reader to \cite[Appendix I]{padic-sigma} and \cite[Ch.~II]{lang}.

The relevance of division polynomials is that they appear in a simple formula for the coordinates of $mQ$ in terms of the coordinates of $Q$. For an integer $m \geq 1$, and a non-torsion point $Q$, we have
\begin{equation}
 x(mQ) = \frac{\theta_m(Q)}{\psi_m(Q)^2},
  \qquad 
  y(mQ) = \frac{\omega_m(Q)}{\psi_m(Q)^3},
\end{equation}
where $\theta_m, \omega_m \in \QQ(E)$ are certain auxiliary functions defined in terms of the $\psi_m$.

The quantities $\psi_m(Q), \theta_m(Q), \omega_m(Q) \in \QQ$ will generally not be integers, due to the coordinates of $Q$ themselves having denominators. It is convenient to introduce a normalising factor that absorbs these denominators. Accordingly we set
\[
\begin{aligned}
 \norm\psi_m(Q) & = \psi_m(Q) d(Q)^{m^2-1}, \\
 \norm\theta_m(Q) & = \theta_m(Q) d(Q)^{2m^2}, \\
 \norm\omega_m(Q) & = \omega_m(Q) d(Q)^{3m^2}.
\end{aligned}
\]
Note that $\norm\psi_m$, $\norm\theta_m$ and $\norm\omega_m$ are defined only on $E(\QQ)$ --- they are most definitely \emph{not} rational functions on $E$. One checks, by examining the degrees of $\psi_m$, $\theta_m$ and $\omega_m$, that $\norm\psi_m(Q)$, $\norm\theta_m(Q)$ and $\norm\omega_m(Q)$ are all integers. Therefore we now have \emph{two} representations of $x(mQ)$ and $y(mQ)$ as ratios of integers,
\begin{equation} \label{eq:two-ratios}
\begin{aligned}
 x(mQ) & = \frac{\alpha(mQ)}{d(mQ)^2} = \frac{\norm\theta_m(Q)}{\norm\psi_m(Q)^2 d(Q)^2}, \\
 y(mQ) & = \frac{\beta(mQ)}{d(mQ)^3} = \frac{\norm\omega_m(Q)}{\norm\psi_m(Q)^3 d(Q)^3}.
\end{aligned}
\end{equation}

The point of \eqref{eq:two-ratios} is twofold. First, Proposition IV.2 of \cite{wuthrich} guarantees that, under the hypothesis of Proposition \ref{prop:division-polys}, the fractions on the right are in fact \emph{reduced} fractions; in other words, that $d(mQ) = \pm \norm\psi_m(Q) d(Q)$. Therefore we may conclude from \eqref{eq:two-ratios} that
\begin{equation} \label{eq:division-poly-equiv}
\begin{split}
 d(mQ) & = \pm \norm\psi_m(Q) d(Q), \\
 \alpha(mQ) & = \norm\theta_m(Q), \\
 \beta(mQ) & = \pm \norm\omega_m(Q),
\end{split}
\end{equation}
where the choices of signs in the first and third equations agree.

Second, we will show that $\norm\psi_m(Q)$, $\norm\theta_m(Q)$ and $\norm\omega_m(Q)$ can be efficiently computed mod $R$ using the usual recursion formulae for the division polynomials. For our application, it is not necessary to compute the division polynomials themselves --- in fact, to do so would completely miss the point, since their degree grows like $m^2$, which is precisely the rate of growth that we are trying to avoid. Rather, we need only compute their \emph{values at $Q$}, and only mod $R$. Fortunately the standard recursive formulae for division polynomials, with minor modifications, are perfectly suited to this task.

To avoid losing $p$-adic precision, we give a version of the formulae that involve no divisions (apart from one division by $2$, which is permitted since we have assumed that $2$ is invertible in $R$). Also, our formulae are tailored to computing the normalised versions $\norm\psi_m(Q)$, $\norm\theta_m(Q)$ and $\norm\omega_m(Q)$  directly, instead of $\psi_m(Q)$, $\theta_m(Q)$ and $\omega_m(Q)$, so that we can work with integral quantities throughout.

In the formulae below, we abbreviate $\alpha(Q)$ by $\alpha$, and similarly with the other variables. We start by defining normalised versions of the coefficients $a_k$ of the elliptic curve, setting
 \[ \norm a_k = d^k a_k, \quad k = 1, 2, 3, 4, 6. \]
We next define variables $\norm b_k$ and $\norm B_k$, which are normalised versions of the $b_k$ and $B_k$ appearing in \cite{padic-sigma}:
\[
\begin{aligned}
 \norm b_2 & = \norm a_1^2 + 4\norm a_2, \\
 \norm b_4 & = \norm a_1 \norm a_3 + 2\norm a_4, \\
 \norm b_6 & = \norm a_3^2 + 4\norm a_6, \\
 \norm b_8 & = \norm a_1^2 \norm a_6 + 4\norm a_2\norm a_6 - \norm a_1 \norm a_3 \norm a_4 + \norm a_2 \norm a_3^2 - \norm a_4^2, \\
 \norm B_4 & = 6\alpha^2 + \norm b_2 \alpha + \norm b_4, \\
 \norm B_6 & = 4\alpha^3 + \norm b_2 \alpha^2 + 2\norm b_4 \alpha + \norm b_6, \\
 \norm B_8 & = 3\alpha^4 + \norm b_2 \alpha^3 + 3\norm b_4 \alpha^2 + 3\norm b_6 \alpha + \norm b_8.
\end{aligned}
\]

Similarly, we need normalised versions of the $g_m$ from \cite{padic-sigma}, defined by
 \[ \norm g_0 = 0, \quad \norm g_1 = 1, \quad \norm g_2 = -1, \quad
   \norm g_3 = \norm B_8, \quad
   \norm g_4 = \norm B_6^2 - \norm B_4 \norm B_8, \]
and then recursively for $m \geq 5$ by
\begin{equation}
\label{eq:recursive}
\begin{split}
\norm g_{2n+1} & =
    \begin{cases}
       \norm B_6^2 \norm g_{n+2} \norm g_n^3 - \norm g_{n-1} \norm g_{n+1}^3, & \text{$n$ even}, \\
       \norm g_{n+2} \norm g_n^3 - \norm B_6^2 \norm g_{n-1} \norm g_{n+1}^3, & \text{$n$ odd}, \\
    \end{cases} \\
 \norm g_{2n} & = \norm g_n(\norm g_{n-2} \norm g_{n+1}^2 - \norm g_{n+2} \norm g_{n-1}^2).
\end{split}
\end{equation}

Finally, the values of $\norm\psi_m$, $\norm\theta_m$ and $\norm\omega_m$ for $m \geq 2$ are given in terms of the $\norm g_m$ by
\begin{equation}
\label{eq:final}
\begin{aligned}
 \norm T & = 2\beta + \norm a_1 \alpha + \norm a_3, \\
 \norm\psi_m & = \norm T^{\sigma(m+1)} \norm g_m, \\
 \norm\theta_m & =
     \alpha \norm\psi_m^2 - \norm\psi_{m+1} \norm\psi_{m-1}, \\
  \norm\omega_m & = \frac{-1}2 \left(
      \norm T^{\sigma(m)}
         (\norm g_{m-2} \norm g_{m+1}^2 - 
          \norm g_{m+2} \norm g_{m-1}^2)
      + \norm \psi_m (\norm a_1 \norm \theta_m + \norm a_3 \norm\psi_m^2)
       \right),
\end{aligned}
\end{equation}
where $\sigma(k)$ is $0$ or $1$ accordingly as $k$ is even or odd.

The algorithm implementing Proposition \ref{prop:division-polys} now runs as follows. All computations are performed mod $R$. We are given as input $\alpha(Q)$, $\beta(Q)$ and $d(Q)$, and the constants $a_k$. From \eqref{eq:division-poly-equiv} it suffices to compute $\norm\psi_m$, $\norm\theta_m$ and $\norm\omega_m$.

We start by computing all of the $\norm a_k$, $\norm b_k$ and $\norm B_k$, and $\norm g_0$ through $\norm g_4$, using the formulae given above. Then, using \eqref{eq:recursive}, recursively compute $\norm g_{m-2}$ through $\norm g_{m+2}$. During this recursive step, it is important to retain the values of $\norm g_j$ as they are computed, since many of them will be reused. Finally, the equations \eqref{eq:final} then determine $\norm\psi_m$, $\norm\theta_m$ and $\norm\omega_m$.

Now we analyse the complexity. Arithmetic operations in $\ZZ/R\ZZ$ may be performed in time $\OO(\log R)$, so the crux of the matter is to show that the recursive formulae \eqref{eq:recursive} are evaluated at most $O(\log m)$ times.

Let $k \geq 4$. To determine $\norm g_j$ for all $j$ in the range $k \leq j \leq k + 7$, using \eqref{eq:recursive} it suffices to know $\norm g_j$ for $(k-3)/2 \leq j \leq (k+11)/2$ if $k$ is odd, or $(k-4)/2 \leq j \leq (k+10)/2$ if $k$ is even. In other words, to determine 8 consecutive values of $\norm g_j$ near $j = k$, it suffices to know 8 consecutive values of $\norm g_j$ near $j = k/2$. Iterating this process, to compute $\norm g_k$ one must evaluate \eqref{eq:recursive} at most $8\log_2(k) = O(\log k)$ times.

\begin{example}
\label{sec:division-polys-example}

Let $P = (5/4, -3/8)$ be a generator of $E(\QQ)$ where $E$ is the elliptic curve $y^2 + y = x^3 + x^2 - 7x + 5$. This curve is `91b1' in Cremona's database, and has conductor $91 = 7 \cdot 13$. One checks directly that $P$ reduces to a non-singular point at the bad primes $7$ and $13$.

We will illustrate Proposition \ref{prop:division-polys} for the case $R = 99$, $m = 101$. The starting point $P$ is specified by
 \[ \alpha = 5, \quad \beta = 96, \quad d = 2. \]
The various constants are initialised as
\[
\begin{aligned}
\norm a_1 & = 0,   & \norm b_2 & = 16,  & \norm B_4 & = 6,  \\
\norm a_2 & = 4,   & \norm b_4 & = 73,  & \norm B_6 & = 4,  \\
\norm a_3 & = 8,   & \norm b_6 & = 57,  & \norm B_8 & = 67,  \\
\norm a_4 & = 86,  & \norm b_8 & = 59,  & \norm T & = 2. \\
\norm a_6 & = 23,
\end{aligned}
\]

We now list the intermediate results in the computation of $\alpha(mP)$, $\beta(mP)$ and $d(mP)$ mod $R$. The required values of $\norm g_m$ are
\[
\begin{aligned}
\norm g_0 & = 0, & \norm g_6 & = 63, & \norm g_{12} & = 0, & \norm g_{23} & = 35, & \norm g_{48} & = 0, & \norm g_{99} & = 49, \\
\norm g_1 & = 1, & \norm g_7 & = 98, & \norm g_{13} & = 64, & \norm g_{24} & = 0, & \norm g_{49} & = 1, & \norm g_{100} & = 19, \\
\norm g_2 & = 98, & \norm g_8 & = 35, & \norm g_{14} & = 71, & \norm g_{25} & = 91, & \norm g_{50} & = 62, & \norm g_{101} & = 82, \\
\norm g_3 & = 67, & \norm g_9 & = 50, & \norm g_{15} & = 4, & \norm g_{26} & = 17, & \norm g_{51} & = 49, & \norm g_{102} & = 72, \\
\norm g_4 & = 10, & \norm g_{10} & = 73, & \norm g_{16} & = 1, & \norm g_{27} & = 67, & \norm g_{52} & = 46, & \norm g_{103} & = 98. \\
\norm g_5 & = 37, & \norm g_{11} & = 98, & \norm g_{22} & = 1, & \norm g_{28} & = 46, & \norm g_{53} & = 1, \\
\end{aligned}
\]
From these it follows that
\[ \norm \psi_{100} = 38, \quad \norm \psi_{101} = 82, \quad \norm \psi_{102} = 45, \]
and
\[ \norm\theta_{101} = 32, \quad \norm\omega_{101} = 4. \]
Therefore we obtain
\[
\begin{aligned}
\alpha(101P) & = 32, \\
\beta(101P) & = \pm 4, \\
d(101P) & = \pm 2 \cdot 82 = \pm 65. \\
\end{aligned}
\]

The result may be verified by using a machine to explicitly compute the coordinates of $101P$; it turns out that the negative sign was the correct one.
\end{example}


\section{A complete example}
\label{sec:example}

We will illustrate the main algorithm by stepping through a computation of the $p$-adic height of a rational point on the curve
 \[ E\colon y^2 + xy = x^3 - 12x + 16, \]
which is `214a1' from Cremona's database. It has conductor $214 = 2 \cdot 107$, and its rank over $\QQ$ is one, with generator
 \[ P = (0, -4). \]
To make life interesting, we will take $p = 43$, which is anomalous for $E$ since $\# E(\FF_{43}) = 43$. We will aim to compute $h_p(P)$ to precision $p^6$, so set $M = 6$. The Tamagawa numbers of $E$ at $2$ and $107$ are respectively $7$ and $1$; therefore we set
\[ \begin{aligned}
  n_1 & = \# E(\FF_{43}) = 43, \\
  n_2 & = \LCM(7, 1) = 7, \\
  n & = \LCM(n_1, n_2) = 301, \\
  m & = n / n_2 = 43, \\
  M' & = M + 2v_{43}(n) = 8.
\end{aligned} \]

To apply Theorem \ref{thm:height} we need $\sigma_p(t)$ mod $I_9$; from Theorem \ref{thm:sigma} this means we require $\EE_2$ modulo $p^6$. To be able to apply Kedlaya's algorithm, it is convenient to change coordinates to put the curve into the Weierstrass form $y^2 = x^3 + ax + b$; for our curve this turns out to be
 \[ E'\colon y^2  = x^3 - \frac{577}{48} x + \frac{14689}{864}. \]
(Note that the equation $E'$ has the same discriminant $\Delta = -13696$ as our original equation, so the value of $\EE_2$ we compute for $E'$ will be the correct value also for $E$. If instead we chose an equation whose discriminant differed by a factor of $u^{12}$, then we would need to adjust the computed $\EE_2$ by $u^2$, since $\EE_2$ is of weight two. In other words, we need to take into account the choice of invariant differential implied by the choice of equation, as $\EE_2(E, \omega)$ is a function of both $E$ and $\omega$.)

We apply Kedlaya's algorithm --- omitting the details --- to find that the Frobenius matrix is
 \[ F = \begin{pmatrix} 4996923274 & 3651910366 \\ 1002107518 & 1324439776 \end{pmatrix} \pmod{p^6}. \]
Then
 \[ F^6 = \begin{pmatrix} 3987851820 & 4837860471 \\ 1528699020 & 2333368599 \end{pmatrix} \pmod{p^6}, \]
so
 \[ \EE_2 = -12 \cdot \frac{4837860471}{2333368599} = 5899790810 \pmod{p^6}. \]
 
Using $\EE_2$ as input, Theorem \ref{thm:sigma} then yields the sigma function,
\begin{multline*}
 \sigma_p(t) = t + (135909305554 + O(p^7)) t^2 + (3933286396 + O(p^6)) t^3 \\
     + (129848206 + O(p^5)) t^4 + (2650487 + O(p^4)) t^5 \\
     + (77893 + O(p^3)) t^6 + (1561 + O(p^2)) t^7 + (8 + O(p)) t^8 + O(t^9).
\end{multline*}
(See Example \ref{ex:sigma} for a more detailed example of computing $\sigma_p(t)$.)

Now we bring into the picture the particular $P$ whose height is sought. We compute
 \[ Q = n_2 P = \left( \frac 34, \frac{-25}8 \right), \]
so that $\alpha(Q) = 3$, $\beta(Q) = -25$, and $d(Q) = 2$. Using Proposition \ref{prop:division-polys}, we then find that, mod $p^8$,
\[ \begin{aligned}
 \alpha(mQ) & = 9491762277279, \\
 \beta(mQ) & = \pm 10171094217691, \\
 d(mQ) & = \pm 3360349669562. \\
\end{aligned} \]
(See Example \ref{sec:division-polys-example} for a more detailed example of this step.) Substituting everything into \eqref{eq:height-2} and \eqref{eq:sigma-subst}, we find that
\[ \begin{aligned}
 h_p(P) & = \frac{2}{301^2} \log_p(1430987165464 + O(43^8)) \\
   & = \frac{2}{7^2 \cdot 43^2} (43 \cdot 44668563676 + O(43^8)) \\
   & = 43^{-1} \cdot 96127622779 + O(43^6),
 \end{aligned} \]
to precision $p^6$ as desired.


\section{Sample computations}
\label{sec:samples}

The author implemented the above algorithms in the computer algebra system SAGE \cite{sage}, building on an implementation of the algorithm of \cite{MST} by Robert Bradshaw, Jennifer Balakrishnan, Liang Xiao, and the author. The code is freely available under a GPL license, and is distributed as a standard component of SAGE (version 2.2 and later). An example session:
\begin{verbatim}
sage: E = EllipticCurve("37a"); E
 Elliptic Curve defined by y^2 + y = x^3 - x over Rational Field
sage: P = E.gens()[0]; P     # a generator of E(Q)
 (0 : 0 : 1)
sage: h = E.padic_height(p=5, prec=5)
sage: h(P)
 4*5 + 3*5^2 + 3*5^3 + 4*5^4 + O(5^5)
\end{verbatim}

The examples below were run on an AMD Opteron running at 1.8 GHz, kindly provided by William Stein, funded by his NSF grant \#0555776. For these examples, we did \emph{not} use the methods of \S\ref{sec:E2-improvement}, opting instead to use the trace and determinant as a correctness check.

\subsection{Large prime case}

We will take the elliptic curve `92b1' from Cremona's database, which has equation $y^2 = x^3 - x + 1$ and conductor $92 = 2^2 \cdot 23$. A generator of the group of rational points is $P = (x, y) = (1, 1)$.

Let $p = 10^{11} + 3$ and $M = 6$. We need to compute $\EE_2$ mod $p^4$. Using the method of \cite{kedlaya-large-p}, SAGE finds that the matrix of Frobenius on the standard basis for Monsky--Washnitzer cohomology is given by
 \[ \begin{matrix}
 \scriptstyle 64304585760876115698175680344198318130013083 &
 \scriptstyle 42503972380936025561124602310186734870477220 \\
\scriptstyle  97128558385368210540568141789457735335273547 & 
\scriptstyle 35695414251123884302364319655812481869943749
\end{matrix} \]
mod $p^4$. The computation time was 42 minutes. As a consistency check, one may verify that the determinant of this matrix is $p \pmod{p^4}$, and that the trace is
 \[ 100000000012000000000540000000010799999956832 \equiv -43249 \pmod{p^4}, \]
which agrees with other point-counting algorithms.

From the matrix SAGE finds immediately that
  \[ \EE_2 = 74470168280485533213508423470741122284560152 + O(p^4), \]
and that the $p$-adic sigma function is
\begin{multline*}
 \sigma_p(t) = t + (69769590353020230550922850977954746761856727 + O(p^4)) t^3 \\ + (5214659177355434704657 + O(p^2)) t^5 + O(t^7).
\end{multline*}

Finally, SAGE uses Theorem \ref{thm:height} to find that $h_p(P)$ is
 \[ p \cdot 9226324270539878944369124959203473806055293044599072658 + O(p^6). \]
This last step is virtually instantaneous. Observe that the original algorithm of \cite{MST} would have required computing the coordinates of the point $nP$ where $n \approx 10^{11}$, which would occupy around $10^{22}$ bits of storage --- not to mention the computation time. Clearly, Proposition \ref{prop:division-polys} is essential for handling such large $p$.

\subsection{High precision case}

Continuing with the same curve, we now consider the case $p = 5$ and $M = 3000$. Therefore we require $\EE_2$ modulo $5^{2998}$. Since $p$ is small relative to $M$, SAGE selects Kedlaya's original algorithm \cite{kedlaya} for the Frobenius matrix computation. The first and last few digits of $\EE_2$ are
 \[ \EE_2 = 3 + 2\cdot5 + 2\cdot5^3 + 3\cdot5^5 + 2\cdot5^7 + \cdots + 3\cdot5^{2995} + 3\cdot5^{2996} + O(5^{2998}). \]
The computation time to obtain $\EE_2$ was 229 seconds.

Obtaining the $p$-adic sigma function took 158 seconds. There are about $3000$ coefficients, each with a similar number of $5$-adic digits, so we dare not write it down here. For such high precision computations, the original algorithm of \cite{MST} would have been utterly impractical; the $M^4$ contribution from the initial power series expansions would multiply the above running time by a factor of perhaps $10^7$.

Finally, computing $h_p(P)$ itself from the sigma function took about 6 seconds. It turns out to be
 \[ h_p(P) =  3\cdot5 + 3\cdot5^2 + 2\cdot5^3 + 5^4 + \cdots + 4\cdot5^{2998} + 2\cdot5^{2999} + O(5^{3000}). \]

\bibliographystyle{amsalpha}
\bibliography{fast-heights}

\providecommand{\bysame}{\leavevmode\hbox to3em{\hrulefill}\thinspace}
\providecommand{\MR}{\relax\ifhmode\unskip\space\fi MR }
\providecommand{\MRhref}[2]{%
  \href{http://www.ams.org/mathscinet-getitem?mr=#1}{#2}
}
\providecommand{\href}[2]{#2}
\begin{thebibliography}{MST06}

\bibitem[Ber]{bernstein}
Daniel Bernstein, \emph{Fast multiplication and its applications},
  http://cr.yp.to/lineartime/multapps-20041007.pdf, retrieved 4th March 2007.

\bibitem[BK78]{brent-kung}
R.~P. Brent and H.~T. Kung, \emph{Fast algorithms for manipulating formal power
  series}, J. Assoc. Comput. Mach. \textbf{25} (1978), no.~4, 581--595.

\bibitem[Bre76]{brent}
Richard~P. Brent, \emph{Multiple-precision zero-finding methods and the
  complexity of elementary function evaluation}, Analytic computational
  complexity (Proc. Sympos., Carnegie-Mellon Univ., Pittsburgh, Pa., 1975),
  Academic Press, New York, 1976, pp.~151--176.

\bibitem[Har07]{kedlaya-large-p}
David Harvey, \emph{Kedlaya's algorithm in larger characteristic}, to appear in
  Int. Math. Res. Not., arXiv preprint \texttt{math.NT/0610973v2}, 2007.

\bibitem[Ked01]{kedlaya}
Kiran~S. Kedlaya, \emph{Counting points on hyperelliptic curves using
  {M}onsky-{W}ashnitzer cohomology}, J. Ramanujan Math. Soc. \textbf{16}
  (2001), no.~4, 323--338.

\bibitem[Lan78]{lang}
Serge Lang, \emph{Elliptic curves: {D}iophantine analysis}, Grundlehren der
  Mathematischen Wissenschaften [Fundamental Principles of Mathematical
  Sciences], vol. 231, Springer-Verlag, Berlin, 1978.

\bibitem[MST06]{MST}
B.~Mazur, W.~Stein, and J.~Tate, \emph{Computation of p-adic heights and log
  convergence}, Documenta Math. (Extra Volume: John H. Coates' Sixtieth
  Birthday) (2006), 577--614.

\bibitem[MT91]{padic-sigma}
B.~Mazur and J.~Tate, \emph{The {$p$}-adic sigma function}, Duke Math. J.
  \textbf{62} (1991), no.~3, 663--688.

\bibitem[Sil92]{silverman}
Joseph~H. Silverman, \emph{The arithmetic of elliptic curves}, Graduate Texts
  in Mathematics, vol. 106, Springer-Verlag, New York, 1992, Corrected reprint
  of the 1986 original.

\bibitem[SJ05]{sage}
William Stein and David Joyner, \emph{Sage: System for algebra and geometry
  experimentation}, Communications in Computer Algebra (ACM SIGSAM Bulletin)
  \textbf{39} (2005), no.~2, 61--64.

\bibitem[Wut04]{wuthrich}
Christian Wuthrich, \emph{The fine selmer group and height pairings}, Ph.D.
  thesis, Cambridge, 2004.

\end{thebibliography}

\end{document}